\documentclass[11pt, oneside]{article}   	
\usepackage{geometry}                		
\usepackage{graphicx}				
				
\usepackage{amssymb}

\title{Galilei and Huygens:  Music and science}

\author{Athanase Papadopoulos\footnote{\emph{Institut de Recherche Mathématique Avancée} and
\emph{Centre de Recherche et d'Expérimentation sur l'Acte Artistique  (ITI CREAA)}. Address: IRMA,
Université de Strasbourg and CNRS,
7 rue René Descartes,
67084 Strasbourg Cedex France;
email: athanase.papadopoulos@math.unistra.fr}}
 \date{}							

\begin{document}
\maketitle

\begin{abstract}

 Vincenzo Galilei and Constantijn Huygens were both  humanists and eminent musicians, the former from the late Renaissance and the latter from the early Modern era. Their respective sons, Galileo and Christiaan, were scientists whose importance cannot be overestimated. 
My aim in this chapter is to set the scene for a parallel presentation of the legacy of the Galilei on the one hand, and  the Huygens on the other. This will give us an opportunity to talk about mathematics, music and acoustics, but also about science in general, at this time of birth of the Modern era. 

The last version of this paper will appear in the book ``Vincenzo Galilei, the Renaissance between music and science", ed. Natacha Fabbri and Ferdinando Abbri, Leo S. Olschki editore,  Firenze, 2025.

\bigskip

\noindent{\bf AMS codes.}  01A45, 00A65,  	85-03.

\bigskip

\noindent{\bf Keywords.}  Christiaan Huygens, Constantijn Huygens, Vincenzo Galilei,  Galileo Galilei, Marin Mersenne, René Descartes, mathematics and music, vibrating strings, falling bodies, curves in the plane, acoustics, musical temperament, history of astronomy, history of music, pendulum.

 \end{abstract}

\bigskip

 \section{Introduction}

Vincenzo Galilei (c. 1520-1591) occupies a special place in the history of Western thought, for his work and his relationship with his contemporaries, and in particular for the education he provided and the impetus he gave to his son Galileo Galilei (1564-1642), the founder of modern physics. In this chapter, I would like to broaden the perspective on Vincenzo Galilei, and at the same time focus on certain scientific and cultural facts that allow us to draw a parallel between, on the one hand, the legacies of Vincenzo Galilei and that of his son Galileo, and on the other, those of Constantijn Huygens (1596-1687) and his son Christiaan (1629-1695).   Addressing the legacy of the Galilei and the Huygens will also give us the opportunity to talk about science and music in general.

Vincenzo Galilei and Constantijn Huygens were excellent musicians, both playing the lute and the keyboard. They were also prolific composers. Huygens' music is still performed today. Vincenzo Galilei was not only an excellent lute player, but he also used this instrument to conduct experiments on musical intervals which supported the theoretical questions he addressed on octave division and the perception of consonance and dissonance. Their respective children, Galileo and Christiaan, developed, each in his own way, his father's ideas on music.  Some authors have rightly insisted on the fact that the young Galileo witnessed his father's experiments on this instrument, and helped him with his mathematical knowledge.\footnote{See e.g. {\sc Stillman Drake},  Music and Philosophy in Early Modern Science, in V. Coelho (ed.), Music and Science in the Age of Galileo, Dordrecht, Kluwer, 1992, pp. 3-16.} We shall also review some of their scientific achievements, in this period of the nascent modern era, presenting these contributions in parallel.  

Mathematics, physics and astronomy are the main fields of science that underwent considerable development in the 16th and 17th centuries.  There is also a field intimately linked to music: acoustics, a theoretical and experimental science that studies sound, its formation and propagation, instruments, musical intervals and several related subjects. Vincenzo Galilei and Constantijn Huygens were both involved in this science. Their children, Galileo and Christiaan, were also engaged in the science of sound, as were other prominent scientists of their time, such as Johannes Kepler (1571-1630),  
René Descartes (1596-1650), Marin Mersenne (1588-1648), John Wallis (1616-1703), Isaac Newton (1643-1727) and several others.  Let me add to this list Blaise Pascal (1623-162), the religious mystic and mathematician who was the contemporary of Christiaan Huygens, one of the main characters of the present chapter: Pascal's first book was on the theory of sound: \emph{Traité des sons}, which he wrote at age 11.
 It is worth noting that the great figures of that epoch who contributed to the study of sound are, broadly speaking, the same as those who made the great mathematical and physical discoveries. In particular, Descartes is regarded as the father of analytical geometry, Newton as one of the two founders of differential calculus --- the other being Gottfried Wilhelm Leibniz (1646-1716) who, incidentally, also took a close interest in the theory of sound ---, Mersenne made substantial contributions to number theory, and one could name several other scientists in this context. In the sections that follow, we shall mention some of Galileo Galilei's and Christiaan Huygens' works on acoustics, and thus touch on the field of music, the favourite activity of their parents.

  \section{Vincenzo Galilei}
  
 Vincenzo Galilei was born in a family of merchants that belonged to the minor Florentine nobility. He grew up in Santa Maria a Monte, a commune in the province of Pisa. He married Giulia Venturi degli Ammannati, who came from a prosperous and cultured background. Vincenzo played the lute, taught music and composed madrigals. He wrote books on theoretical music. He formed his own ideas about singing in general, preferring monodic singing and advocating a sung recitative in the style of ancient Greek theater --- although little was known about the latter. This style was to give rise a few decades later to Italian opera.   Like most musicians of his time, Vincenzo Galilei was, in the pure Greek tradition, interested in the mathematical side of music. He was a member of the \emph{Camerata}, an academy founded by Giovanni de' Bardi (1534-1612), Count of Vernio, a Florentine nobleman, writer, composer, art critic and intimate of the Medici. Bardi was about twenty-four years Galilei's junior, and took him under his protection.  
 
 Of particular note among Galieli's treatises is his \emph{Dialogo della musica antica et moderna}.\footnote{\label{f:Dialogo} {\sc Vincenzo Galilei}, Dialogo della musica antica e della moderna, Florence, G. Marescotti, 1581, Transl. by Claude V. Palisca,  Dialogue on ancient and modern music,  In: Humanism in Renaissance musical thought, New Haven and London, Yale University Press, 1985.}   He is credited with being the first to observe that to obtain the octave of a note emitted by a stretched string, the tension of the string must be multiplied by four and not by two, contradicting the Pythagorean creed that the ratio of integers 1/2 is universally associated with the octave. Vincenzo Galilei describes such an experiment in his \emph{Discorso intorno alla diversita delle forme del Diapason}, written around 1589, a treatise in which he examines the different ways of producing the octave (called, in ancient Greek terminology, the diapason, hence the title of the book).\footnote{\label{f:Diapason} See  {\sc  Vincenzo Galilei}, Discorso intorno alla diversita delle forme del Diapason, Transl. Claude V. Palisca: A special discourse concerning the diversity of the ratios of the diapason, In: The Florentine Camerata: documentary studies and translations, Yale University Press, New Haven, 1988.} His experiment consisted in attaching a system of weights to one end of a vibrating string and observing that in order to obtain the octave while keeping the length constant, the weights must be quadrupled, not doubled. Similarly, to obtain the fifth, one must multiply the weights by 4/9, and not by 2/3.  Vincenzo Galilei's discovery is a consequence of the general law of string vibration stated by Marin Mersenne a few decades later, which says that the frequency of vibrations produced by a stretched string is inversely proportional to the square of the string's tension, and not to the tension itself. Mersenne's law is sometimes considered to be the first known non-linear law of physics.  We shall mention Mersenne several times in the sequel.
     
     Vincenzo Galilei conducted his own experiments. In this context, as in others, he taught his son Galileo the 
importance of experimentation, and of not relying blindly on established theories, no matter how old they may be.      Through his experiments, we get a glimpse of Galilei's ``technician" side, a technician in the original sense of the Greek word \emph{technê}, meaning both ``art" and ``craft".  The following statement, extracted from the same treatise,\footnote{Cf. p. 181 of the \emph{Discorso}, Footnote \ref{f:Diapason}.} has been considered as 	a indication of Galilei's adherence to the ``experimental method" that was to take over science a few decades later: ``There are few things that cannot be weighed, counted or measured."
     
     Galilei's experiments led him to distrust dogmas, or at least not to blindly accept them, even those of well-established authorities, as he wrote in the \emph{Dialogo della musica antica et della moderna}:\footnote{Cf. p. 269 of the English translation of the \emph{Dialogo}, Footnote \ref{f:Dialogo}.}

\begin{quote}\small
In connection with [the theories of Pythagoras] I wish to point out two false opinions of which men have been persuaded by various  writings, and which I myself shared until I ascertained the truth by means of experiment, the teacher of all things.
\end{quote}

Talking about the Pythagoreans, it is worth recalling here that they admitted only rational numbers (ratios of integral numbers) in their musical theories, which rules out all considerations of tempered scales, since these involve taking $n$-th roots of integers, which are generally irrational.\footnote{It may be useful to recall that the Pythagoreans developed a science and philosophy in which integral numbers and their ratios (the ``rational" numbers) are at the basis of the relationships between all kinds of beings, and not only musical. Integral numbers were at the basis of all the Pythagorean theories: from moral and political to cosmic and supercosmic.}  

%
%
%

The same posture of Vincenzo Galilei towards science appears in the following passage, extracted from his \emph{Discorso intomo all'uso delle dissonanze}:\footnote{Transl. in   {\sc Claude V. Palisca}, \emph{Was Galileo's Father an
Experimental Scientist?} In:  Music and Science in the Age of Galileo,
Dordrecht/Boston/London, Kluwer Academic Publishers, 1992, p. 143-151 (the sentence is on p. 143).} 

\begin{quote}\small

Before your Lordship begins to untie the knot of the proposed
questions [concerning the nature of the diatonic practiced today], I
wish in those things which sensation can reach that authority always
be set aside (as Aristotle says in the Eighth Book of the
Physics), and with it the tainted reason that contradicts any [sense] perception at all of truth. For it seems to me that those who, for the sake of proving some conclusion of theirs, want us to believe them purely on the basis of authority without adducing any further
arguments are doing something ridiculous, not to say (with the
Philosopher\footnote{That is, with Aristotle.}) acting like silly fools.


\end{quote}

As for Vincenzo's son Galileo, there is no need to recall how he stood up against dogma; we know that he paid for it with his freedom, having been confined to house arrest for the rest of his life.

I shall now turn to the two Huygens, and I would like to start with three quotations, one from the father and the other two from the son, which illustrate the fact that both were, like Vincenzo and Galileo Galilei, wary of the indiscriminate use of rules.  The first quotation is extracted from a letter from Constantijn Huygens to Henry du Mont, dated October 10, 1658. He writes:\footnote{{\sc Jacob Adolf Worp}  De briefwisseling van Constantijn Huygens (1608-1687), 's Gravenhage,  Martinus Nuhoff, 1911-1917. See Part V, p. 311, letter 5591.}
%
 
 \begin{quote}\small I am not unfamiliar with the rules but I find only little consistency, and many contradictions between the authors   [...]
 
 \end{quote}
%
%

Similarly, his son Christiaan wrote the following sentence in the margin of his manuscript \emph{Cycle harmonique par la division de l'octave en 31 dièses, intervalles égaux}  (1691) (\emph{\OE uvres Complètes}, t. XX, p. 17):\footnote{All the references to Huygens' Complete Works are to {\sc Christiaan Huygens}, \emph{\OE uvres complètes}, 22 vol., Société hollandaise
des sciences, La Haye, 1888-1950.} 
 ``Non audio qui allegant authoritatem."\footnote{Portefeuille Musica, f. 16-19.  \emph{\OE uvres complètes},  t. XX, p. 155-164, with French translation in  p. 327-336 of  the volume {\sc Renzo Caddeo, Xavier Hascher, Franck Jedrzejewski, Athanase Papadopoulos},   Christiaan Huygens: \'Ecrits sur la musique et le son, \oe uvres traduites du latin, avec commentaires math\'ematiques, historiques musicaux, Hermann, Coll. Th\'eorie de la musique, Paris,  2021. In the rest of this this chapter, we shall refer to this volume as \emph{\'Ecrits sur la musique}.}
(I do not listen to those who, to impose themselves, appeal to authority).\footnote{The authority to which Huygens refers in this marginal note is Marin Mersenne.}

In the same vein, let me quote the following sentence by the same Huygens, reprinted in t. XXIII of his \emph{\OE uvres Complètes}, p. 342): ``Doubt is painful to the mind. This is why everyone willingly agrees with those who claim to have found certainty."

    Now we turn in some detail to the second musician we are particularly interested in in this chapter, Constantijn Huygens.
  
  \section{Constantijn Huygens}

 Constantijn Huygens belonged to the Dutch aristocracy. We learn from his biographers that in his childhood he was taught to read and write several languages, together with mathematics and philosophy. This was not unusual for the education of gifted highborn boys. In addition, he received a consequent musical education, first from his father, and then, with famous teachers.  He became a fine player of lute, viola da gamba, theorbo, spinet and organ.
In his autobiography, he writes:\footnote{See {\sc Astrid de Jager}, \emph{Constantijn Huygens' Passion: Some Thoughts About the Pathodia sacra et profana},  Tijdschrift van de Koninklijke Vereniging voor Nederlandse Muziekgeschiedenis, Deel 57, no. 1 (2007), p. 29-41, and Huygens' biography in: Jacob Adolf Worp: Fragment eener autobiografie van Constantijn Huygens in Bijdr. en Meded. van het Hist. Gen. XVIII (1897) p. 1-121; New edition: The Correspondence of Constantijn Huygens [De briefwisseling van Constantijn Huygens 1608–1687], ed. by Ineke Huysman, Huygens ING; with metadata from the Worp edition taken originally from Epistolarium, CKCC project, Huygens ING, The Hague.} ``I have always dedicated myself with constant diligence to the most important part of music, the polyphony. I even dare to say, I have come close to the Italians, whom one can take as a model."
  
%
%
%
%

 Huygens became a diplomat by profession, evolving among noblemen, ambassadors and ministers.  In 1625, he became secretary to Frederick-Henri of Nassau, Prince of Orange, who at the time was an ally of France.  He remained loyal to this princely family throughout his life, representing it in Paris and Versailles on numerous diplomatic missions that also took him  to England and other European countries.  Besides playing several instruments, he was a composer and music theorist. The editors of  Constantijn Huygens' \emph{Correspondence and musical works} report that he placed music above everything else, and that he never traveled without his lute.\footnote{\label{f:Jonc} See the Preface by   {\sc Willem Jozef Andries Jonckbloet},  In: {\sc Jan Pieter Nicolaas Land}  (ed.), Musique et musiciens au XVIIe siècle : Correspondance et \oe uvres
musicales de Constantijn Huygens, Société pour l'histoire musicale des
Pays-Bas, E. J. Brill, Leyde, 1882.}
 
Constantijn Huygens was also a prolific poet, writing in Latin, Flemish, Italian, English, German and French. He translated books written in all these languages, as well as others; he used to say that translating was his way of perfecting his language skills. His voluminous correspondence is written almost entirely in French. His poems and musical compositions remain a significant part of the Dutch cultural heritage.
The editors of Constantijn Huygens' Correspondence and Musical Works estimate the number of his musical compositions at over 800.\footnote{p. XXVI of the Preface by {\sc Jonckbloet}, Footnote \ref{f:Jonc}.} English translations of Constantijn Huygens' poems are available in print.\footnote{Cf. A Selection of the Poems of Sir Constantijn Huygens (1596-1687), In: Amsterdam Studies in the Dutch Golden Age, Amsterdam University Press, several editions.} There are also several recordings of his compositions.

Constantijn Huygens was also a Latinist, with a keen interest in philosophy and mathematics.  Condorcet, at the beginning of his \emph{\'Eloge} of Constantijn's son Christiaan, speaking of his father, writes:\footnote{{\sc Marie-Jean-Antoine-Nicolas de Caritat de Condorcet}, \emph{\'Eloge
d'Huyghens}, in \OE uvres de Condorcet, t. II, publiées par A. Condorcet
O'Connor et M. F. Arago, Firmin Didot, Paris, 1847, p. 54-72.} ``In those times when the reborn enlightenment aroused even more gratitude than envy, it was not considered that a taste for letters and philosophy was in a magistrate a useless or dangerous distraction."

Huygens was close to René Descartes, who had gone into exile in Holland to avoid censorship and he defended Descartes when the latter was accused of impiety. Condorcet, in the same \emph{\'Eloge}, writes: ``We will not praise Constantijn Huygens for having dared to cultivate Latin poetry, and to study the physics and geometry of Descartes; but we will praise him for having had, although a man of position, the courage to be publicly the friend of this philosopher, falsely accused of impiety, and persecuted by a dangerous enemy." 
Huygens and Descartes maintained a lively correspondence, discussing mathematics, dioptrics, mechanical arts and lenses. The two men shared ideas on matter, gravity, the causes of motion and many other topics. As friends, they also discussed their daily preoccupations and they talked about the people they met. Pages 736 to 824 of Volume III of Descartes' \emph{\OE uvres} in  Adam and Tannery's edition is comprised exclusively of letters exchanged between Descartes and Constantijn Huygens   between January 3, 1640 and June 26, 1643.\footnote{ {\sc René Descartes}, \OE uvres
publiées par Charles Adam et Paul Tannery, Paris, Cerf, 
1897-1913.
}  

Descartes visited Huygens regularly. He was impressed by the talents of his friend's son Christiaan, whose inventiveness and passion for mathematics were evident from an early age. Meanwhile, Constantijn kept his friend informed of his son's progress.  Descartes kept a close eye on young Christian's development. He also appreciated his interest in music theory, a subject that had occupied him personally since his youth and which he continued to cultivate through his epistolary exchanges. A case in point is a letter to Constantijn Huygens in 1646\footnote{\emph{\OE uvres Complètes}, t. XXII, letter no. 25.} in which Descartes, after informing the latter of a correspondence he had exchanged with the composer and music theorist Joan Albert Ban,\footnote{Joan Albert Ban (or Joannes Albertus Bannius) (1597-1644) was a Catholic priest, composer and music theorist. He was in contact with Constantijn Huygens, Mersenne and Descartes, and is the author of a \emph{Dissertatio epistolica de musica natura} (1637).} wrote to him about his son Christiaan, who was 17 at the time: 

\begin{quote}\small
If your second son wishes to practice in this matter [music], he may have the opportunity to do so by correcting Bannius and myself and showing that neither of us has understood anything about it. For our reasons being neither mathematical nor physical but only moral, as I said there, it is easy to find others that are contrary to them. If he writes anything on this subject, I shall be very glad to see it.
\end{quote}

Marin Mersenne also regularly corresponded with Constantijn Huygens. He called him ``the father or patron of Music".\footnote{\emph{\OE uvres Complètes} t. I, letter no. 27, quoted below.}     The correspondence between the two men, and later between Mersenne and Christiaan Huygens,  includes many questions related to Galileo Galiei's work, as we shall recall later in this chapter. The Minime was well acquainted with the writings of the two Galileis and in his correspondence he acted as a link between them and the Huygens.  As a matter of fact, Mersenne published a free translation of and commentary on Galileo's \emph{Discorsi e dimostrazioni matematiche intorno a due nuove scienze attenenti alla meccanica e i movimenti locali} in which he commented extensively on the latter's experiments on sound.\footnote{Cf. {\sc Marin Mersenne} (editor and translator) Les nouvelles pensées de Galilée, mathématicien et ingénieur du duc de Florence,  où il est traité de la proportion des mouvements naturels et violents, et de tout ce qu'il y a de plus subtil dans les méchaniques et dans la physique. H. Guenon, Paris, 1639, 256 pages. New edition:  Les Nouvelles Pensées de Galilée, édition critique avec
Introduction et notes par Pierre Costabel et Michel-Pierre Lerner, Avant-propos de
Bernard Rochot, 2  volumes, Coll. L'histoire des sciences, Textes et \'Etudes, Paris, Vrin, 1973.} In his correspondence with Constantin and Christiaan Huygens, Mersenne addressed  several times the question of falling bodies and other questions that Galileo Galilei had dealt with at length. We shall discuss these questions in the sequel.

Among the many subjects covered in the correspondence between Mersenne and his friend Constantijn Huygens is, of course, music. Their exchange on this subject included both technical and practical matters: the two men used to share their opinions on musicians, instruments, music teachers of their time, etc., and also on how to improve music in all aspects. 
 
In a letter dated January 12, 1647,\footnote{\emph{\OE uvres Complètes} t. I, letter no. 27.} Mersenne discussed with Constantijn Huygens the question of harmonic tones (without, however, using the word ``harmonic''), explaining that this phenomenon constituted the greatest difficulty he had ever encountered in music: 
\begin{quote}\small
What and when and why the voice which is among the lowest that one can take, or several voices taken together, make besides their own tone another tone which is above at the twelfth or double fifth. And this also happens to the thick strings of a touched viol. For if you stand in great silence in  your cabinet and with great attention you touch [the viol's string]  gently or faintly [\ldots]  you will always hear a sound [making an] echo accompanying the natural one, at the twelfth above, and often another one at the seventeenth.
 [\ldots] Certainly since you can rightly be called the father or patron of music, this difficulty deserves that you should not leave it to the world, and that, for yourself, it should not be demonstrated by anyone [else than you], [in a manner] so clear that anyone can see the reason for it."\footnote{The last sentence reads in French: ``Certes puisqu'à bon droit l'on peut vous appeler père ou patron de la musique, cette difficulté mérite que vous ne la laissiez pas au monde, et que, pour vous-même, elle ne soit démontrée par qui que ce soit, si clair que chacun en voie la raison."}
 \end{quote}

 In another letter to his friend, dated March 17, 1648 (t. I, letter no. 4), Mersenne discusses questions relating to the optimal shape of the lute (should parabolic or hyperbolic lines be used? etc.). He mentions ``a lute or 15-row theorbo where the quarter tones are marked". On April 6, 1648,\footnote{\emph{\OE uvres Complètes},  t. II, letter no. 47a.}   Constantijn Huygens announces to Mersenne that an English lord had ``brought him from France a fine old grand lute from Bologna, the best he has ever touched". He asks him in turn what form of lute he thought would produce the most beautiful resonance, and why. He tells him, regarding this question: ``I know something by experience, which rarely deceives me, but my Archimedes will be in charge of reasoning on the theory".\footnote{Archimedes was the nickname which Constantine Huygens used for his son Christiaan.}

 Constantijn Huygens was a true art lover. The editors of his son Christiaan's \emph{\OE uvres Complètes}  report that Constantijn was one of the first (perhaps the first) to understand Rembrandt's genius.\footnote{\emph{\OE uvres Complètes} t. XXII, p. 405.}  He was responsible for commissioning paintings for the House of Orange, and he bought several Rembrandts for them. 
   
 We turn now more specifically to the two offsprings, Galileo Galilei and Christiaan Huygens.
 
 \section{Galileo Galilei and Christiaan Huygens, I}
 
    Galileo Galilei and Christiaan Huygens were both raised in environments conducive to their intellectual development. As children, they both had the opportunity to meet musicians and other artists who passed by the family home. Vincenzo Galilei wanted his son to become a physician, while Constantijn Huygens thought his son would follow in his footsteps into a career in diplomacy. Neither son followed in his father's path: Galileo became a professor of mathematics, then of mechanics and astronomy, while Christiaan lived mainly on his father's annuities, then on those of Louis XIV, who took him under his protection, and later on from some of his inventions' income, principally watches and clocks, for which he had obtained exclusive exploitation privileges from the King of France.

In their approach to science, Christiaan Huygens and Galileo Galilei shared a geometrical bent when it came to tackling mathematical and physical problems, in the pure tradition of the ancient Greeks and particularly Archimedes. This applies to fields such as optics, hydrostatics and, more generally, the science of equilibrium of bodies, to problems  of areas and volumes, and many others. The two men never met. When Galileo died in 1642, exiled in his villa in Arcetri, Huygens was just thirteen years old.

From a certain point onwards, Mersenne, who already had a sustained correspondence with Constantijn Huygens, began, at the latter's request, to write to his son. Indeed, in a letter to Mersenne dated December 23, 1646,\footnote{\emph{\OE uvres Complètes}, t. I, letter no. 23a.} Huygens writes: ``From now on, I will only serve as an address to the literary trade among you two''. His son, whom he called ``my little mathematician", was 17 at the time. A few days later (letter dated January 14, 1647,\footnote{\emph{\OE uvres Complètes}, t. I, letter no. 27a.}Huygens writes to Mersenne: ``You have too high an opinion of my Archimedes, but I know that he will never deny the one you must have of him." The letter is a reply to one from Mersenne dated January 3, 1647,\footnote{\emph{\OE uvres Complètes}, t. I, letter no. 24.} in which the latter comments on Christiaan's ideas on falling bodies and in which he tells his friend that he believes that if his son continues to work in the way he started, he will one day surpass Archimedes. The problem of falling bodies was of particular interest to Galileo Galilei, and we shall comment on this below.

 Christiaan Huygens was primarily a mathematician. He was a founding member of the Académie Royale des Sciences, and probably the most notable Parisian mathematician of his time. Gottfried Wilhelm Leibniz (1646-1716), during a stay he made in Paris from 1672-1676, was his private pupil.\footnote{Leibniz arrived to Paris in 1672 as a diplomat, sent by the Elector of Mainz on a mission to Louis XIV.  He was 26 years old and had just invented his computing machine. The secret reason for which Leibniz accepted this mission was that he wanted to make contact with a number of scientists, including Christiaan Huygens, who was seventeen years his senior, who had been living in Paris for six years. Leibniz' plan succeeded, and he studied geometry and physics with Huygens throughout his stay in Paris.} The two men always maintained a rich epistolary relationship.  In their correspondence, Leibniz often presented his mathematical works to Huygens, asking for his opinion.
 
%
%
%

%
 In the following sections, I shall discuss various scientific topics of interest to both the Galilei and Huygens families. For chronological reasons, I will consider questions concerning vibrating strings, falling bodies, the pendulum, curves and astronomy in that order; these are all questions in which Huygens was interested from his youth. I will then turn to the problem of musical temperament, a question which was dear to the two Galilei and the two Huygens, and on which Christiaan Huygens did substantial work later in his lifetime.

I begin with the question of string vibration. For Galileo Galilei, this was a continuation of the work started by his father. For Christiaan Huygens, the work was also a tribute to his own father, whose favorite activity was music.

\section{Galileo Galilei and Christiaan Huygens, II: The vibrating strings}

%
%
%

 Christiaan Huygens was well aware of Galileo Galilei's work on the theory of sound, and he mentions it in his piece \emph{Le Son},\footnote{\emph{\OE uvres Complètes}, t. XIX, p. 353 ff.} Part II of this piece is titled: ``Rapports des longueurs des cordes consonantes suivant Pythagore, et rapports des nombres de leurs vibrations suivant Galilée et d'autres
savants."\footnote{The reference to Galielo's work on this subject is the \emph{Dialogo Secundo} of Galileo Galilei's \emph{Discorsi e dimostrazioni matematiche intorno a due nuove scienze attenenti alla meccanica e i movimenti locali},   Leida, 1638, Ludovico Elzeviro.}
 
 Mersenne, in a letter to Christiaan Huygens dated November 16, 1646,\footnote{\emph{\OE uvres Complètes}, t. I, letter no. 17.} addresses a related question. He writes:

 \begin{quote}\small
 I forgot to ask [your father] if you can play  
the lute; if so, I beg you to see if you would like this
beautiful problem, why the string $AB$, taken at your wish, for example
chanterelle,\footnote{The chanterelle is the thinnest and sharpest string in a stringed instrument.} attached firmly at $A$ and making some
sound, must be stretched at $B$ four times as much as before to reach the
octave, whereas it only needs to be shortened by half at $C$ to bring it up to
the octave. I understand that your foundation of mechanics teaches that to make a movement twice as fast,
perhaps a quadruple force is required; you will make the beginning of the
proof, which will be invaluable to me from your hand.
 \end{quote}

One may recall here that Mersenne, in his \emph{Harmonie universelle} published ten years earlier (1636),\footnote{ {\sc Marin Mersenne},  Harmonie universelle, contenant la théorie et la pratique
de la musique, Sébastien Cramoisy, Paris, 1636.} gave a formula which (in modern terms) says that whereas the frequency of a sound produced by a vibrating string is inversely proportional to its length,  the same frequency is proportional to the square root of the string tension. In fact, in his letter to Huygens, Mersenne seeks to understand the \emph{reason} for the need of taking the square root of the tension.  
%
 This question is also linked to the discovery made by Vincenzo Galilei which we already mentioned, namely, that in order to make a string sound at the upper octave, the tension of this string must be divided by the factor 1/4 and not by 1/2. This kind of problem concerning musical intervals continued to occupy Huygens for the rest of his life. Huygens replied to Mersenne a few days later (the letter is dated November 1646):\footnote{\emph{\OE uvres Complètes}  t. I, letter no 20.}
\begin{quote}\small 
Concerning the musical problem which you suggested to me, \emph{amplius deliberandum censeo}:\footnote{\emph{Ego amplius deliberandum } is a famous sentence by Terence (2nd c. BC), which was also used by others after him, meaning that the question is unclear and needs more thought, information or deliberation.} having found it in your book physicomat., I have often speculated on it; but the solution is very difficult as far as I can see, and it is necessarily so, because otherwise it would not have been ignored by so many brave minds today.
\end{quote}
 
It is important to recall here that Mersenne performed a great number of experiments, which he describes in his 
\emph{Harmonicorum libri, in quibus agitur de sonorum natura, causis, \& effectibus}.\footnote{{\sc Marin Mersenne}, Harmonicorum libri, in quibus agitur de sonorum natura, causis, \& effectibus:
de Consonantiis, Dissonantiis, Rationibus, Generibus, Modis, Cantibus, Compositione,
orbisque totius Harmonicis Instrumentis, Paris, Lutetiae Parisiorum, Guillaume Baudry, 1636.} Among his propositions, let me mention the one saying that the tone produced by a vibrating string does not depend on the amplitude of the vibration (this is Proposition XXIX, p. 24). We may find this natural, but this is not a priori obvious.

In his letters to Christiaan Huygens, Mersenne asks questions, and at the same time he shares his thoughts on the problem he was addressing. In his letter dated January 8, 1647,\footnote{\emph{\OE uvres Complètes},  t. I, letter no. 25.} he tackles once again the question of the influence of the string tension on the musical intervals produced. He introduces it as follows: ``[Your father] is very fond of music, and I hope you will take pleasure, for his sake, in solving a harmonic difficulty that you yourself will try out on your lute [. . . ]". He then submits to him the same question we mentioned before (the dependence of tone on tension), in greater detail, explaining how the tension of a string can be doubled by using weights.
In the same letter, Mersenne writes: ``Regarding what I told you about the quadruple tension of the harmonic string to bring it up to the octave, I remember that I dealt with this in the \emph{Ballistics},\footnote{See \emph{Ballistica, et acontismologia}, part of the  \emph{Cogita physico-mathematica}, Paris, 1644, p. 1-140 (special pagination).} Proposition 36, the reading of which will perhaps bring the true demonstration to your mind."

Several years later, in 1673, Christiaan Huygens discovered a profound relationship between the problem of the vibrating string and that of the pendulum, on which he was working extensively. Thanks to this relationship, he was able to establish the mathematical laws that are behind Mersenne's formula on the frequency of a vibrating string. We refer the interested reader to Huygens' article \emph{Découverte de la théorie générale de l'isochronisme des vibrations}.\footnote{\label{f:ecrits} \emph{\OE uvres Complètes}  t. XVIII, p. 489-495, also reproduced in Chapter 16 of \emph{\'Ecrits sur la musique}.}

 Another problem that Christiaan Huygens tackled and which had attracted Galileo's interest before him is that of falling bodies. The problem concerns a search for the law of acceleration of these bodies, but also, and more fundamentally, it is about understanding why objects we drop start to fall. Again, this is not an obvious fact. This will be discussed in the next section.
 
  \section{Galileo Galilei and Christiaan Huygens, III: The falling bodies}
  
We start by recalling a few facts which should help understanding the problem of falling bodies.

During the Renaissance, measurements were made on objects dropped from the top of a tower, and it was noticed --- a fact which is not obvious a priori --- that heavy or very heavy objects dropped simultaneously reach the ground almost at the same time. It was also known that the speed of a falling body is not constant, but undergoes an acceleration, and it was thought that this acceleration takes place ``according to the sequence of natural numbers",\footnote{This was Galileo's terminology, see also Footnote \ref{f:Euclid}.} a sentence which meant that if we divide the duration of the fall of a body into equal periods of time, and if during the first lapse of time the body travels one unit of length, it travels two units in the second lapse of time, three in the third, and so on. In other words, it was thought that velocity was proportional to the time covered, or, equivalently, and using again modern language, that the acceleration of free-falling bodies is constant. Laws of falling bodies are already contained in Aristotelian physics.\footnote{Aristotle, in his \emph{Physics}, Book IV §3, states a law according to which (expressed in modern terms) the distance covered in a time t by a moving point subjected to a force F is equal to the product Ft. We know that this law is false, but at least it has the merit of existing. It is usually considered as the oldest known expression of a law of motion.} Galileo got the idea of experimenting on falling bodies using inclined planes on which he dropped balls, postulating that the law governing falling bodies is the same for an object falling on a vertical plane as on an inclined plane; the advantage of working on an inclined plane is that the object falls more slowly, and thus lends itself better to measurement. In this context, it is important to recall that Galileo did not have a clock with which he could measure time. To study falling bodies, he was led to devise ingenious devices,  such as a water clock, that is, he measured the elapsed time by evaluating the weight of water flowing from a thin pipe to the bottom of a large water jug. Such experiments are described in the \emph{Discorsi e dimostrazioni matematiche intorno a due nuove scienze}.\footnote{Third day of the \emph{Discorsi }, see p. 108. of Vol. III of the \emph{Opere},  {\sc Galileo Galilei}, Le Opere, New edition, G. Barbèra,
Florence, 1966.}

%
%
%

We now turn to Christiaan Huygens, who was interested in the same question from an early age. In a letter dated September 3, 1646\footnote{\emph{\OE uvres Complètes}, t. I, letter no 11.} (he was 17 years old), he wrote to his brother (also called) Constantijn:
\begin{quote}\small
AB is an altitude from which a weight C is dropped. I demonstrate that in the first [interval of] time of its fall it passes a spacing like here CD, in the second [interval of] time equal to the first, 3 such spacings, and comes as far as E, in the third [interval of] time 5 spacings, in the fourth 7, and so will continue to make greater progress each time, always adding to the last spacing twice the first one. ... In addition to this, I have demonstrated that if it is thrown in any direction, it describes a parabola; of all this and an infinite number of other things that depend on it, I never knew the demonstration until I invented it myself. 

\end{quote}
%

Huygen's solution was the same as Galileo Galilei's.
 Mersenne knew that this problem was difficult. In a letter dated October 13, 1646,\footnote{\emph{\OE uvres Complètes}, t. I, letter no. 13a.} he wrote to Christiaan Huygens about this problem, declaring in particular: ``As I greatly honor [your father], and as I believe it will please him in talking with you about your propositions of which you say you have the demonstration, I will only tell you about the last one, of which I do not believe you have the demonstration if I do not see it."  
He then developed his ideas about the falling bodies, without reaching any definite conclusion. At the same time he talks to his young correspondent about ballistics, and the parabolic trajectory of objects. This is a subject to which Galileo had   made significant contributions. Indeed, Galileo was not only interested in the vertical fall of bodies, but also in the trajectory of those to which an initial (non-vertical) impetus is given. He discovered that this trajectory (unless the impetus is directed vertically) is always parabolic, a result that formed the basis of what was to become theoretical ballistics, a topic which was already of interest to Constantijn Huygens (the father), before his son took up the subject in his turn. 
 In the same letter, Mersenne writes to Christiaan Huygens that, for the trajectory to remain parabolic, ``the impetus imparted to the projectile must never cease; but qualities that are easily imparted, such as impetus, are just as easily and quickly lost, \emph{Violentum non durabile}''. To conclude on this subject, he says: 
 \begin{quote}\small

 Nevertheless, if notwithstanding this consideration you believe that your demonstration is still valid, I will be pleased if you communicate it to me; and then, if you have got it right, I will tell you to show this violence [force] so that you can then determine the locus  by which your effect must strike with more violence.  I would add that the principles Galileo has taken from everything he has said about motion are hardly firm, and that although in small heights the proportions follow fairly closely, in large they are almost always lacking.  [\ldots]

 \end{quote}

On November 16, 1646,\footnote{\emph{\OE uvres Complètes}, t. I, letter no. 17.} Mersenne writes again to Christiaan Huygens: ``I assure you that I have so greatly admired the beauty of your demonstration of the falling [bodies] that I believe Galileo would have been delighted to have you as a guarantor of his opinion." In the letter to Christiaan Huygens dated January 8, 1647,  which we have already quoted in connection with the question of string tension,\footnote{\emph{\OE uvres Complètes}, t. I, letter no 25.} Mersenne once again addresses the problem of falling bodies in connection with his recent work \emph{Novarum observationum physico-mathematicarum}\footnote{{\sc Marin Mersenne}, Novarum observationum physico-mathematicarum, Paris, Antoine Bertier, 1647.} that he intends to send him soon after: 
\begin{quote}\small
I assure myself that all this will be well worthy of your consideration, which I hope you will share with me, namely in what it impugns\footnote{This means that in his new work, Mersenne shows that this proportion is not correct.} the proportion of acceleration by Galileo's numbers 1, 3, 5, 7, etc. which is also [yours], and in that it contains many things, which [contradict] certain principles of M. Descartes.\footnote{\label{f:Euclid}  It is important to remember that Galileo did not have the tools of the differential calculus which was developed a few decades later by mathematicians such as Leibniz, Newton, Huygens, Jacob and Johann Bernoulli, and the Marquis de L'Hôpital. In particular, he was not familiar with the notion of derivative to express velocity, and even less with that of second derivative to express acceleration. To formulate his mathematical results, he used the geometrical language of the Greeks, that of Archimedes and Apollonius, often expressing these results with graphical representations, speaking of ``deviation of the tangent at different points of a curve" or using similar expressions.
In the passage in hand, Galileo is relying on the theory of proportions  set out in Euclid's \emph{Elements}, with its limitations. In particular, according to this theory, one is not allowed to divide by each other two quantities that are not of the same nature, for example, a distance and a time; it was therefore impossible for Galileo to define speed as a distance divided by a time, or to write that a distance is proportional to the square of a time. He wrote, for example, that the acceleration of a body in free fall is ``following the sequence of odd numbers: 1, 3, 5, 7 . . .",  meaning that if at time 1 the distance covered is 1, at time 2 it will be 1 + 3 = 4, at time 3 it will be 4 + 5 = 9, at time 4 it will be 9 + 7 = 16, and so on. In modern terms, this means that the distance covered is proportional to the square of the time.}
\end{quote}

Beyond the law of falling bodies, Huygens wanted to understand \emph{why} objects that are dropped start to fall. Even if the phenomenon is not surprising, because we are used to it, the question remains legitimate.  Aristotle wrote extensively on this subject. As with everything he studied, he sought the causes. The explanation he gave was that every object tends towards its natural place --- heavy objects towards the earth, fire towards the sky. Even today, despite all the ancient and modern theories available on this subject, the question of ``why" and ``how" attraction works remains open, or at least very difficult to explain.

 It is natural to see that Pierre-Simon de Laplace (1749-1827), in his \emph{Exposition du système du monde}, retains only Galileo's name for the problem of falling bodies. He writes:\footnote{{\sc Pierre-Simon de Laplace},   Exposition du système du monde. Paris, Imprimerie du Cercle-Social, rue du Théâtre-Français, n° 4. L'An IV de la République française (1796), t. I, p. 238.} ``Galileo laid the first foundations of the science of motion, with his beautiful discoveries on the fall of bodies. Geometers, following in the footsteps of this great man, have at last reduced the whole of mechanics to general formulae that leave nothing to be desired but the perfection of analysis."\footnote{``Galilée jeta les premiers fondements de la science du mouvement, par ses belles découvertes sur la chute des corps. Les géomètres, en marchant sur les traces de ce grand homme, ont enfin réduit la mécanique entière à des formules générales qui ne laissent plus à désirer que la perfection de l'analyse."} But at the same time, he attributes to Huygens the discovery of the relation between the problem of falling bodies and that of the motion of the pendulum  (\emph{Op. Cit.}, t. II, p. 260): `` This transition from the oscillating motion, whose duration can be observed with great precision, to the rectilinear motion of the falling bodies, is an ingenious remark for which we are still indebted to Huygens. " \footnote{``[\ldots] Ce passage du mouvement d’oscillation dont on peut observer avec une grande précision la durée, au mouvement rectiligne des graves, est une remarque ingénieuse dont on est encore redevable à Huygens."} Thus, we shall pass now to the pendulum, whose motion is also part of the field of kinematics, which was at the heart of the work of Galileo Galilei and Christiaan Huygens. The study of the pendulum is closely related to the concept of a clock. Galileo spent the last years of his life thinking of the conception of such an instrument. In building a pendulum clock, Huygens realized one of Galileo's dreams. Let us see this in some detail.

\section{Galileo Galilei and Christiaan Huygens, IV: The pendulum} 

Galileo's name is also associated with the theory of the pendulum. He is considered to be the first to advocate the study of the isochronism of the pendulum, i.e., the regularity of its oscillations, within a rigorous physico-mathematical framework, even though this regularity had certainly been observed long before him. He is said to have had the idea of isochronism in 1581 (he was 17), while he was a medical student: Observing the sway of a candlestick,\footnote{According to a legend, the swing was that of the chandelier in Pisa Cathedral; cf. {\sc Alexandre Koyré}, Galilée et l'expérience de Pise : à propos d'une légende, In:  Études d'histoire de la pensée scientifique, Gallimard, Paris, 1965, p. 213-223.}  he was puzzled by its regularity and he tried to confirm his intuition by counting the period of the oscillation,  comparing it with his own pulse rate.\footnote{It is worth keeping in mind here that, as noted above, in those days, there were no watches or clocks to measure time. Much later, Galileo wrote about this phenomenon in his \emph{Dialogo intorno ai due massimi sistemi del mondo, Tolemaico e Copernicano}, published in 1632. This is Galileo's main writing in which he defended the theory of heliocentrism, a theory which had been condemned in 1615 by the Roman Inquisition. This was of course the book that triggered Galileo's trial, which took place in the year following its publication.}

In 1658, the same year Huygens published the first version of his \emph{Horologium oscillatorium}\footnote{\emph{\OE uvres Complètes}, t. XVII.} and 16 years after Galileo's death, Prince Leopold de' Medici asked Vincenzo Viviani, a disciple of Galileo and later of Torricelli and who held the title of first engineer to the Grand Duke of Tuscany, to write a memoir in which he would formally clarify Galileo's contribution to the pendulum clock.  Four years before, Viviani had published a biography of Galileo,\footnote{Cf.  {\sc Vincenzo Viviani}, Racconto istorico della vita del Sig.r Galileo Galilei di Vincentio Viviani, In: 
A. Favaro, Edizione Nazionale delle Opere di Galileo Galilei, Firenze, Barbera Editore, XIX, App. III, 1968, p. 602 ff. This is the first biography of Galileo. In this work, which predates Huygens' \emph{Horologium} by four years, Viviani makes no mention of any idea of Galilei of using a pendulum to build clocks. Viviani only mentions the laws of the pendulum and how it had been used in medicine to measure pulse beats, and in astronomy to measure the duration of celestial phenomena.} in which he reports on the pendulum and its use in astronomy and medicine (for measuring pulses). But Viviani did not mention there any work on the pendulum clock. In his report to the prince, Viviani claimed that Galileo had imagined such a pendulum, but had never built one. Galileo's son, Vincenzio Galilei,\footnote{Vincenzio Galilei (1619-1649) was an illegitimate son of Galileo Galilei and  Marina Gamba,  with whom he had two other children, likewise illegitimate.} tried to construct one,  with the help of a Florentine craftsman, Domenico Balestri, but the project was halted after Balestri died. Viviani nevertheless mentions the name of a watchmaker working for Grand Duke Ferdinand II, who built a pendulum clock a few years later, following an idea different from that of Vincenzio Galilei. At the same time, he acknowledges that Huygens was the first to build a pendulum clock. 
Volume III of Christiaan Huygens' \emph{\OE uvres Complètes} contains a document written by Viviani addressed in the form of a letter to Leopoldo de Medicis in which he reports on this question.\footnote{Letter no. 673b., August 20, 1659; Appendix I to no. 673a.  In the same volume of the \emph{\OE uvres Complètes}, the editors mention that Girolamo Tiraboschi, in his \emph{Storia della Letteratura Italiana}, Roma, 1785,  t. VII, p. 155, quotes the \emph{Novelle Fiorentine} of 1774 (No. 10, p. 150), in which it is said that ``Senator Nelli owns a story of the pendulum clock, written in 1659 by Vincenzo Viviani, according to which Galileo Galilei only imagined such a clock in 1640, but did not build it, that his  son Vincenzio tried to build one with the help of Domenico Balestri, a Florentine craftsman, but that, surprised by death, he was unable to see it completed."}

 Ismail Boulliau,\footnote{Ismail Boulliau (1605-1694) was a mathematician, astronomer, translator of ancient Greek texts, librarian and world traveler.  He was a close friend of Pierre Gassendi, Christiaan Huygens, Marin Mersenne and Blaise Pascal, and was one of the first foreign associates of the Royal Society.} in a letter dated January 9, 1660,\footnote{\emph{\OE uvres Complètes}, t. III, p. 8, Letter no. 707.} sent Christiaan Huygens a drawing he had received from Florence, depicting a pendulum clock begun by Galileo Galilei.
Huygens replied to Boulliau in a letter dated January 22, 1660 (letter no. 711. p. 12 of the same volume): 

\begin{quote}\small

You have given me great pleasure by sending me the design of the clock that Galileo had begun. I see that the pendulum is as good in it as in mine; but not applied in the same way; for firstly he has substituted a much more awkward invention instead of using the wheel, which is called ``de rencontre". Secondly, he did not suspend the pendulum from a net or small ribbon, but [he did it] in such a way that all its weight rests on the axis on which it moves, which is undoubtedly the main reason why his method has not been successful enough; for I know from experience that the movement becomes much more difficult, and the clock prone to stopping. Even though Galileo had the same thought as me about the use of the pendulum, this is rather to my advantage than otherwise, because I have done what he did not know how to do, and yet I had no clue from him or anyone else in the world that this invention was possible.

\end{quote}

%
The question of Galileo's possible discovery of the pendulum clock was regularly raised.\footnote{See the account in \emph{\'Ecrits sur la musique}, p. 122 ff.}
 In any case, it is certain that Galileo carried out a considerable number of experiments on the pendulum. Not only did he discover its isochronism for low-amplitude oscillations, but he also came up with the idea of using the pendulum's swing as a means of measuring time. He needed urgently such a device for his astronomical observations, for which he used to appoint a person for counting the number of beats. 

To end this section, let me mention an excerpt from a letter addressed   by Christiaan Huygens to Pierre Petit\footnote{Pierre Petit  (1598-1677)  was an engineer, advisor and geographer, then intendant general of fortifications  to Louis XIV. He was associated with Pascal, Huygens and other scientists.} on November 1st, 1658:\footnote{\emph{\OE uvres Complètes}  t. II, letter no. No 546.} ``As for determining the length of the pendulum for each clock, I would have done well to explain in general what proportion the vibrations of various pendulums have to do with their length, but all that concerns these proportions having been fully dealt with by Galileo, Father Mersenne and many other authors, I did not feel it necessary to include it in this description."

 In conclusion, it should be clear that before Huygens, Galileo had the intuition that the pendulum could be used to make clocks; nevertheless, he was unable to realize this project, even though the desire to have an accurate clock haunted him for a long time. 

\section{Galileo Galilei and Christiaan Huygens, V: Curves} 

In this section, I would like to discuss questions related to mathematical curves, in particular two geometrically defined ones, which have applications in art and  science, and which attracted the attention our two scientists, Galileo Galilei and Christiaan Huygens. The first curve goes by many names, catenary, hanging chain and others.\footnote{The curve is referred to in both Mersenne's and Huygens' letters as ``chaîne pendue" or ``corde pendante", that is, ``hanging chain" or ``handing rope", but there are also other names: chaînette, alysoid, etc.} it is a curve on which both Galileo and Huygens worked, and where the latter corrected the work of the former. The second curve is called \emph{cycloid}, and is central to the theory of pendulum clocks. Both curves have also been used in architecture, the first in the construction of large domes, and the second in that of bridges.

We start again with the correspondence between Huygens and Mersenne.

In a letter to Mersenne, dated October 28, 1646,\footnote{\emph{\OE uvres Complètes}, t. I, letter no. 14.} Christiaan Huygens raises another question that Galileo had already addressed, namely the form of the so-called \emph{catenary}, i.e., the chain hanging from its two ends when it subjected only to its own weight.  In his \emph{Discorsi e dimostrazioni matematiche intorno a due nuove scienze}, published in 1638, a work considered one of the founding texts of modern physics, Galileo erroneously asserted that this curve is a piece of parabola. In his letter to Mersenne, Huygens points out Galileo's error. He announces to his correspondent that he intends to send him, in a forthcoming letter, a demonstration of the fact that the curve is not an arc of parabola. At the same time, Huygens writes that he will show Mersenne a calculation of the pressure that must be applied to a ``mathematical or gravity-free string" to make it describe a parabola. Huygens  also foresaw that this curve could not be described by a simple equation.

It might be useful to recall here that mathematicians since Antiquity have studied a large number of curves. The Greeks have extensively  studied the conics (ellipse, parabola, hyperbola), but also curves known as ``mechanical curves", defined as solutions to some problems on geometric loci or some construction problems and requiring mechanical instruments to be constructed. These include the Archimedean spiral,\footnote{This is a curve used in squaring the circle, in the trisection of an angle, and, in the 20th century, in the design of gramophone records.}   the cissoid,\footnote{A degree-three curve, used in the duplication of the cube, studied by Diocles in the 3rd-2nd century BC and by Christiaan Huygens and Pierre de Fermat in the 17th century.} the quadratrix,\footnote{This curve was introduced by Hippias around 420 BC, and it was used in the trisection of an angle.} the conchoid,\footnote{A generic name for curves obtained as curves at a certain finite distance from other curves. The most famous conchoid is attributed to  Nicomedes, 3d c. BC, and it is used in the problem of the trisection of an angle.} and there are many others. Because of their esthetic appeal, all these curves have been used in modern art. 

The catenoid does not seem to have been considered in Greek antiquity. In the analytic geometry that had recently been founded at the time of Christiaan Huygens, the curves that were studied were essentially defined by equations, and the catenoid does not belong to this class. It was natural that this exceptional curve would attract Christiaan Huygens' interest, given his curiosity and unconventional spirit. Huygens mentions his result on the catenoid at the end of a letter to Mersenne dated October 28, 1646:\footnote{\emph{\OE uvres Complètes}, t. I, letter no. 14.} ``I will finish for fear to detain you for too long, and I will send you in another letter the demonstration that a hanging rope or chain does not make a parabola, and what the pressure must be on a mathematical or gravity-free rope to make [a parabola]; [a result] of which I also found the demonstration, not long ago."

Mersenne replied on November 16 of the same year,\footnote{\emph{\OE uvres Complètes},  t. I, letter no 17.} saying that he was looking forward to reading Huygens' demonstration, and that what he claims would indeed contradict what Galileo had believed. He added that he would be delighted to understand, as the young Huygens thought he knew, what pressure was needed to be exerted on the chain to make it take the form of a parabola, and also, he said, what pressure needed to be exerted on the string to make it take the shape of a hyperbola or ellipse. ``If you can do this", he writes, ``you will surpass yourself''. Several epistolary exchanges took place between Mersenne and Huygens on this topic.  The latter's letters contain mathematical propositions and demonstrations of the facts he claimed. The interested reader will find these letters in Volume I of the \emph{\OE uvres Complètes} of 
Christiaan Huygens.

To conclude this discussion on the hanging chain, let me mention that this curve, in its inverted form, was used in a number of dome constructions, including that of Santa Maria del Fiore Cathedral in Florence, the intermediate dome of the Panthéon in Paris, and the Chernobyl Arch, inaugurated in 2016, which covers the cracked former sarcophagus. The reason is more than aesthetic: the catenary is an arc whose weight ensures its stability. It plays a central role in the transmission of forces, yet another intuition of artists that has outstripped the formal proof of mathematicians.

  We continue now with another curve which attracted the attention of our two scientists, the cycloid. This curve is obtained as the trajectory of a point located on the circumference of a wheel as the latter moves by rolling along a straight line. This curve satisfies a beautiful property known as tautochronism, which says that if we make a material point slide without friction along this curve, the point will always take the same amount of time to reach the bottom of the curve, regardless of its initial height. This is a fundamental discovery that Christiaan Huygens made while he was working on the pendulum clock.

 Blaise Pascal (1623-1662) formulated several problems concerning the cycloid: equation, arclength, area of the region under this curve, volume of the body generated by its rotation, and so on. Such questions could not fail to interest Christiaan Huygens, as they concerned a mechanically defined curve, a subject he was passionate about. Volume II of his \emph{\OE uvres Complètes} contains notes by Pascal, extracts from his correspondence concerning this curve as well as exchanges on the same subject between Pascal, Ismail Boulliau and other mathematicians. According to C. de Waard the editor of Marin Mersenne's \emph{Correspondence}, the problem of studying the cycloïd was first proposed by Mersenne around the year 1615.\footnote{See  {\sc Cornelis de Waard}, \emph{Une lettre inédite de Roberval du 6 janvier 1637 contenant le premier
énoncé de la cycloïde}, Bulletin des Sciences Math.,  Série II, T. XLV, p. 220-224, Paris, Gauthier-Villars, 1921.}

 Huygens discovered the tautochronism of the cycloid in 1658 and realized that he could use it in his work on the pendulum clock. Indeed, he saw  that if the pendulum was made to swing along such a curve, then the period of oscillation becomes independent of the pendulum length.  With such a device, it is as if the length of the pendulum constantly changes, while the period of oscillation remains the same.  He discovered another important mathematical property of the cycloid, namely, if the pendulum is made to swing between two equal arcs of a cycloid, then its end, instead of describing a circular arc,   describes a cycloidal arc. In mathematical terms, this property means that the involute of a cycloid arc is also a cycloid.

 Galileo studied the cycloid before Huygens, and he liked this curve. He used it for drawing bridge arches. He did not foresee the properties that Huygens found. 
In a letter to Bonaventura Cavalieri dated February 24, 1640, he writes, regarding this curve:\footnote{Letter no. 3972 of {\sc Galileo Galilei}, In: Le Opere, New edition, G. Barbèra,
Florence, 1966.}``It has been more than fifty years since it occurred to me to describe this arched curve, admiring its graceful curvature to adapt it to the arches of a bridge. I have made several attempts to demonstrate some property of this curve and of the space enclosed between it and its chord, and at first I thought that this space could be three times the size of the circle that describes it; but this was not the case, although the difference is not great."\footnote{``Quella linea arcuata sono più di cinquant'anni che mi venne in mente il descriverla, e l'ammirai per una curvità graziosissima per adattarla agli archi d'un ponte. Feci sopra di essa, e sopra lo spazio da lei e dalla sua corda compreso, diversi tentativi per dimostrarne qualche passione, e parvemi da principio che tale spazio potesse esser triplo del cerchio che lo descrive; ma non fu così, benchè la differenza non sia molta."} On the work of Galileo Galilei on the cycloid, we refer the interested reader to Laloubère's edition.\footnote{{\sc Antoine de Laloubère},   De Cycloïde Galilaei et Torricellii Propositiones viginti Autore Lalovera Societatis
Jesu, Toulouse, 21 juillet 1658, Bibliothèque nationale de France, Réserve des livres rares, RÉS-V-857.}

The cycloid, as a curve, had already been studied in Antiquity and had attracted the interest of Renaissance mathematicians. Variations of it are extensively used in modern art.\footnote{The number of pieces of modern art titled ``cycloid" is almost boundless.}

It remains for us to say a few words about Christiaan Huygens the astronomer, as the heir of Galileo Galilei.

\section{Galileo Galilei and Christiaan Huygens, VI: Saturn}

Christiaan Huygens' main astronomical work is his \emph{Systema Saturnium},\footnote{\emph{\OE uvres Complètes}  t. XV.} a treatise published in 1659, in which he reports on his discovery of Saturn's moon and its ring. In 1657, he began observing Saturn with a telescope that was far more powerful than Galileo's: its magnification factor was 100, whereas Galileo's was only 30. 

Huygens dedicated his \emph{Systema Saturnium}, to Prince Leopoldo de Medici. For Huygens, this dedication was another way of declaring himself in the line of Galileo, who had worked in Florence several decades earlier, under the protection of the same princely family, and who, in January 1610, had discovered four satellites of Jupiter, which he also dedicated to the Medici family, naming them Medicea Sidera (Medicean stars) in honour of Cosimo II, Grand Duke of Tuscany, and his three brothers, Francesco, Carlo and Lorenzo.\footnote{Today, these four satellites are known simply as the Galilean satellites.} Galileo's discovery once again went against the Aristotelian belief that all the planets and other celestial objects revolve around the Earth, providing further indication that the Earth is not the absolute center of the astronomical system, and that other planets may have their own satellites like the Earth has.
Our two scientists, Galileo and Huygens, used telescopes they had designed themselves.  Huygens polished his own lenses; in fact, it was one of his favorite occupations.\footnote{This is another story, but it was through his lens-polishing work that Huygens came into contact with Spinoza, who himself made his living from this activity. Talking about lenses, let me mention that Huygens was the inventor of the device known as the \emph{magic lantern}, which used a series of lenses to project images painted on glass onto a wall or screen.}. Galileo Galilei was able to observe Saturn's ring, but he did not realize that it was a ring; he considered it a mysterious object (he called it "mysterious appendages"; in fact, depending on its orientation, the ring appeared in different shapes). It was Huygens, with his powerful telescope and systematic observations, who realized that it was a ring (in fact, a collection of rings).
 It is good to read an excerpt of the beautifully written dedication, addressed by Christiaan Huygens to Prince Leopoldo de Medici:
\begin{quote}\small
In this pamphlet, I study objects that are very far away in celestial space, objects that are outside the field of human observation, unless they benefit from the help of applied science. I am sure that many people will say that I have taken too much trouble to examine things which, according to general opinion, are of little concern to us, while there are still so many of those situated here below, close to us, that deserve to be studied. But those who speak this way seem to notice too little how much the investigation of celestial things is superior to any other study, and how grandiose is the fact itself that our contemplation extends to parts of nature placed at such great distances; which, although appearing obscure and small, are nevertheless in reality brilliant and very large.  
For if we were to consider that these objects have little relevance to us because they are so far away, we would surely be unworthy of a mind endowed with the reason by which we easily ascend through the immensity of celestial spaces, unworthy also of that marvelous and never sufficiently praised instrument invented to extend vision, an instrument by means of which we also reach the region of the stars with the sense of sight itself. It is by taking advantage of this invention that I have now penetrated further than anyone before into the remote domain of Saturn, that I have reached so far that of this immense road only a hundredth part remains: if I had somehow been able to cross this last part, how many, and what news, good gods, would I have to tell!

\end{quote}

Finally, let me quote the beginning of the \emph{Systema Saturnium}\footnote{The text is reproduced, in Latin and in French, in volume XV of the \emph{\OE uvres Complètes}, t. XV, p. 214 ff.} in which Huygens refers to Galileo's observations: ``When Galileo used optical glasses to contemplate the celestial bodies [...] and that, before anybody else, he had made known to mortals the 
remarkable planetary phenomena, the most astonishing  
of his observations seem to me those he published
concerning the planet Saturn."

In the last section, I will consider a question of music theory to which I alluded above, namely, that of temperament.

\section{Music and temperament}\label{s:music}

    Volume XX of Christiaan Huygens' \emph{\OE uvres Complètes} is titled \emph{Musique et mathématique}.  Pages 1 to 174 of this volume are dedicated to music.  It is not for mere convenience that the editors of the \emph{\OE uvres Complètes} have included both subjects in the same volume. Huygens' work on music is based on mathematics, in line with most of his eminent predecessors'  works on the subject. In this respect, let me quote from a letter Christiaan Huygens sent  R. Moray\footnote{Robert Moray (c. 1608-1673) was a Scottish scientist and  one of the founders of the Royal Society.} dated August 1st, 1661:\footnote{\emph{\OE uvres Complètes}, t. III, p. 307-308.}
  \begin{quote}\small
  I spent a few days studying music, and the division of the monochord, to which I happily applied algebra. I also found that logarithms are of great use there, and from that point I began to consider these marvelous numbers and admire the industry and patience of those who gave them to us.
    \end{quote}
     
 Let me mention that Christiaan Huygens' works on the theory of sound are not contained in the volume on \emph{Musique et mathématique}; they are reproduced and commented on in other volumes of the \emph{\OE uvres Complètes}, e.g. Volume XIX, whose general title is \emph{Mécanique théorique et physique, 1666-1695}, which contains his notes on his experiments on the human voice, on vibrating strings on the propagation of sound and on other related topics, and Volume XVIII whose title is \emph{L'horloge à pendule ou à balancier de 1666 à 1695} in which part of his work on harmonic vibration is reproduced.
    
     Among the works included in Volume XX of the \emph{\OE uvres Complètes} are essays on singing and others that contain valuable insight into the works of ancient and modern authors and the various modes, tetrachords and temperaments they used. Two pieces are of particular interest to us here, each of which is a collection of notes devoted to temperament; these are \emph{La division du monocorde}\footnote{\emph{\OE uvres Complètes},  Vol. XX, p. 41-60. For the French word ``monocorde", Huygens writes ``monochorde", which is closer to the Greek.}   and \emph{Le nouveau cycle harmonique}.\footnote{\emph{\OE uvres Complètes},  Vol. XX, p. 139-174.}  The former, originally written in Latin, is translated into French in chapter 8 of the volume \emph{\'Ecrits sur la musique} dedicated to Christiaan Huygens' work on music, with notes and comments by the editors.\footnote{See Footnote \ref {f:ecrits}.} This memoir consists of detailed computations, based on the use of logarithms, of a large number of temperaments and octave divisions. Some of these temperaments are divisions with quarter and seventh comma.  Huygens' main aim is to find a temperament that best approximates the accuracy of the intervals of fifths, fourths, thirds and sixths.      
     
          The second collection of pieces, \emph{Nouveau cycle harmonique}, contains the foundations of Huygens' theory on the division of the octave into 31 equal intervals.  Huygens recommends this temperament for both instruments and voice. He shows that the differences between the intervals of this temperament and those of the mesotonic\footnote{The mesotonic (or meantone) temperament is a temperament whose construction is based on the desire to have as many perfect thirds as possible.} are indistinguishable. The pieces that make up this last essay contain highly detailed calculations of musical intervals. This work, which is probably Christiaan Huygens' best-known contribution to music theory, had a non-negligible impact in twentieth century music.\footnote{The impact of Huygens' work on the division of the octave into 31 equal intervals is discussed discussed at some length in chapters 17 and 18 of \emph{\'Ecrits sur la musique}, which also covers some developments of the 31-degree temperament theory both before and after Huygens. Particular mention is made of instruments such as Nicola Vicentino's archicembalo and arciorgano, Vito Trasuntino's clavemusicum omnitonum, Giovanni Battista Doni's cembalo pentarmonico, Francesco Nigetti's cembalo onnicordo, as well as  keyboard musical instruments built in the 20th century, capable of producing Huygens' 31 sounds, reviewed by Adriaan Fokker (1887-1972), the Dutch physicist and composer who wrote music according to Huygens' temperament of 31 equal intervals. Fokker also discussed the relationship between Huygens' work on music and Leonhard Euler's work on the use of the harmonic seventh and the number 7, and his  classification of genres.   The same chapter includes an annotated description of several twentieth- and twenty-first-century musical compositions based on these
31 sounds, works that demonstrate the direct impact of Christiaan Huygens' theoretical work on our modern music.}
 
 I would like to mention the piece \emph{Notes se rapportant à des écrits de musicologues anciens}.\footnote{\emph{\OE uvres Complètes}, t. XX, p. 89-103.}  The full article in French translation can be found in \emph{\'Ecrits sur la musique}, p. 203-226. In this article, Huygens provides a critical review of some of the musical works of ancient authors, such as Eratosthenes' \emph{Division of the Monochord}, Marcus Meibom's tables, the theories of Aristoxenus, Euclid, Gaudentius, Ptolemy, Aristides Quintilian, Wallis and many others. He mentions Gioseffo Zarlino (1517-1590) on several occasions, in relation with his views on the Ancients. I will focus on the latter, as Vincenzo Galilei was a vocal critic of his theories. 
     In fact, I will quote Huygens on Zarlino's work from another piece, which follows the one just mentioned and which is titled \emph{Notes se rapportant à des écrits sur la musique d'auteurs modernes}. In this context, the ``modern" era is the one that begins with Francisco Salinas (1513-1590) and Zarlino. This article is longer than the first, and is presented, in French, on pages 227-265 of \emph{\'Ecrits sur la musique}.\footnote{The original is mixed, French and Latin.} In this article, Huygens mentions the division of the octave into 12 equal semitones, in a paragraph that begins with: ``Aristoxenus divided the octave into 12 equal semitones, which Vincenzo Galilei maintains is the best division. Mersenne uses it for the lute."    In several places in this article, Huygens criticizes Zarlino's theories. In particular, he writes:\footnote{See p. 223 of \emph{\'Ecrits sur la musique}.}
   \begin{quote}\small
   
  Zarlino, in \emph{Ragionamento} 4, Proposition 1, says that [the notion of] temperament was discovered by someone else --- without knowing who it was --- and by chance. Similarly, in Book IV, chap. 12 of the \emph{Sopplimenti musicali} (34), he highly praises this invention. In the \emph{Istitutioni harmoniche}, book II, chap. 42, he describes a mediocre temperament in which major and minor thirds are equidistant from perfect pitch, and fifths and fourths are 2/7 of a comma apart. And in chap. 43, he tries to show that there is no such thing as an acceptable temperament, and disapproves of the one in which 1/2 comma is subtracted from the major tone and the same is added to the minor. This temperament is, however, true and excellent.  At the time, however, this had not yet been fully understood. Indeed, Zarlino mentions that [only the] fifths and fourths remain just there. [\ldots]
     \end{quote}
      
      Huygens' criticism refers to certain passages in Zarlino's \emph{Dimostrationi harmoniche} and others in his \emph{Istitutioni harmoniche}.  In particular, the last sentence of this passage refers to p. 128 of the \emph{Istitutioni harmoniche},\footnote{{\sc Gioseffo Zarlino}, Le Istitutioni harmoniche, nelle quali, oltra le materie appartenenti alla musica, si trovano dichiarati molti luoghi di poeti, d'historici, \& di filosofi,  si come nel leggerle si potrà  chiaramente vedere, Francesco de' Franceschi, Venise, 1558; reprint, Broude Brothers, New York, 1965.} where Zarlino alludes to a temperament in which the syntonic comma would be distributed between the major and minor tones, shortening one by half a comma and lengthening the other by another half a comma, and in which all other intervals --- and therefore fourths and fifths --- are left in their natural proportion. However, in the true mesotonic temperament the fifth is lowered by a quarter of a comma and the fourth raised by as much. We have commented on all these criticisms in the volume \emph{\'Ecrits sur la musique}.\footnote{See Footnote \ref{f:ecrits}.}       A little further on, Huygens writes of Zarlino: ``He makes a big deal of what he knows of geometry, which only concerns proportions, as do all the others, like the doctors in music Boetius, Glarean and Salinas'', and again, that Zarlino ``has muddled everything by mixing ancient music and its terms with modern, as in naming the tones by the names of the Greeks hypate hypaton, trite diezeugmenon, etc., and always talking about tetrachords''.   This is followed by more technical criticisms of Zarlino's proposed interval arithmetic, in which Huygens points out the latter's errors.          
        
     Huygens summarized in a synthetic way the 31-part equal temperament he recommends in his well-known \emph{Lettre à Basnage de Beauval touchant le cycle harmonique (known as Novus Cyclus harmonicus}, written in 1691.\footnote{The letter, which was published in the \emph{Histoire des Ouvrages des sçavans} of which Henri Basnage de Beauval was the director (Rotterdam, October 1691, p. 78-88), is reproduced on p. 169-174 of t. X of the \emph{\OE uvres Complètes}. } The letter is accompanied by calculations and tables which show that the difference between the intervals of this division and those of the fifth (corresponding to the proportion 3/2), the fourth (4/3), the major third (5/4), the minor third (6/5), the major sixth (5/3) and the minor sixth(8/5), are so small that it is ``completely impossible for the most delicate ear" to perceive them.
In the same letter, Huygens notes that Salinas mentions this division of the octave into 31 equal parts, ``but only to condemn it". Similarly, he says, Mersenne rejects it. In conclusion, he claims to have ``demonstrated the excellence of this division by the principles of geometry, and to have upheld it against the unfair judgment pronounced by these two famous writers."

\medskip
 
  Let me conclude this chapter by a general thought on science.

\section{In guise of a conclusion}

We usually talk about the \emph{Scientific revolution} to designate the period of Western science around the end of the sixteenth century and the first decades of the seventeenth (roughly speaking corresponding to the period when Galileo Galilei was most active). 
The term ``revolution" used in this formula is meant to suggest a break with the ``old" and a radically new approach to the principles and processes of science, with a stress on experimentation.\footnote{I would like to quote here Arkady Plotnitsky, who highlights the fact that the stress on mathematics is just as important, who prefers the term ``transformation" to ``revolution" and from whom I borrow a statement by Heidegger: ``Modern science is experimental because of its mathematical project". For the source of Heidegger's quote, and also some reason for this view, see the article {\sc Arkady Plotnitsky}, “All the resources of pure mathematics”: Mathematical physics and mathematics as physics, In: Essays on Topology -- Dedicated to Valentin Poénaru, ed. L. Funar and A. Papadopoulos, Springer Cham, 2025, p. 505-576.} In this sense, in the word ``science", music is arguably included, at least for what concerns acoustics. Let me insist on the fact that even with a ``revolution" taking place,  the discoveries themselves remained very slow; they often involved the works of several generations of scientists. This is one of the conclusions that can be drawn from the account I have given here.

 \bigskip

  \bigskip
  
  \bigskip

\noindent {\bf Acknowledgements} This work is supported by the lnterdisciplinary Thematic lnstitute CREAA, as part of the ITI 2021-2028 program of the Université de Strasbourg, the CNRS, and the Inserm (funded by IdEx Unistra ANR-10-IDEX-0002, and by SFRI-STRAT'US ANR-20-SFRI-0012 under the French Investments for the Future Program).  The author would like to thank  Arkady Plotnitsky who read a first version and made useful comments, and the Erwin Schr\"odinger Institute (Vienna) for its hospitality during the time where this chapter was written.  
   \end{document}